\documentclass[11pt,a4paper]{article}
\usepackage[latin1]{inputenc}
\usepackage{amsmath}
\usepackage{amsfonts}
\usepackage{amssymb}
\usepackage{makeidx}
\usepackage{graphicx}

\newtheorem{thm}{Theorem}
\newtheorem{prop}{Proposition}

\makeindex

\begin{document}

\begin{center}
\textbf{\Huge{Locally conformally cocalibrated $G_2$-structures}} \\
\vspace{15mm}
\Large{Arezoo Zohrabi} \\
\vspace{1cm}
\Large{Dipartamento di matematica, Universit\`{a} degli studenti di Torino, Italy}\\
\end{center}
 \vspace{1cm}

 \begin{center} \textbf{Abstract}\end{center}
  We study the condition in which $G_2$-structures are introduced by a non closed four-form, although they are satisfying locally conformal conditions.All solutions are found in the case when the Lee form of $G_2$-structures is non-zero and $\mathfrak{g}$ introduces seven-dimensional Lie algebras, The main results are given in preposition1 and theorem1.\vspace{2cm}

\tableofcontents

 \newpage \begin{center}\section*{Introduction} \end{center}
  A $G_2$-structure $(M^7,\varphi)$ on a seven-dimensional manifold $M^7$ can be characterized by the existence of globally 3-form $\varphi$ called fundamental 3-form. the classes of $G_2$-structures can be described in the terms of the exterior derivatives of fundamental 3-form $\varphi$ and the 4-form $ \phi=*\varphi$ , where $*$ is the Hodge operator defined from metric and the derivative of $V_{g_{\varphi}}$, in this paper we focus our attention on the class of locally conformally cocalibrated $G_2$-structures (usually denoted by L.C.CC $G_2$-structures), which are characterized by the condition $d\phi=-\theta \wedge \phi$ for a closed non vanishing 1-form $\theta$ also known as the Lee form of the $G_2$-structures.
  Let M be connected n-dimensional manifold with a Riemannian metric on M, Let $\bigtriangledown$ be the Levi-Civita connection of G and for $p,q \in M$ joint by smooth path $\gamma$ then parallel transport along $\gamma$ using $\bigtriangledown$  define an isometry between the tangent space of $T_{p} M $ , $T_{q} M$.\\
  A seven-dimensional connected,oriented Riemannian manifold $M^7$ with the Holonomy  
    \footnote{The holonomy group of metric g is the group of isometries of $T_{p} M $ generated by parallel transport around closed loop at the point p, i.e. $Hol(g)\leq O(n)$}contained in $G_2$ is characterized by the existence of a $G_2$ reduction of its orthogonal frame bundle with the restriction of Cartan form of the Levi-Civita connection, alternatively can be seen that $G_2 \subset SO(7)$ is defined as the stabilizer under the action of $GL(7,\mathbb{R})$ of stable 3-form $\varphi~{'} $ or equivalently 4-form $\phi ^{'} $ which are related $\phi ^{'} = *\varphi~{'}$ \\  
  A seven-dimensional nilmanifold is a compact manifold which corresponding lie group G acts transitively on and in homogeneous case $G/\varGamma$ where $\varGamma$ is a discrete subgroup whit the property for $a.g \in G/\varGamma \Longrightarrow b(a.g)=(ba).g \in G/\varGamma$ for $a,b \in \varGamma$ and $g \in G$, equivalently Lie group G can define a nilmanifold if and only if its lie algebra $\mathfrak{g}$ admits a basis $\{e_1,e_2,...,e_7\}$ and invariant forms on a homogeneous case are determined by forms on $\mathfrak{g}$, so we can consider seven-dimensional nilpotent Lie algebras consequently.\\
  $G_2$ is a compact lie group as a closed subgroup of the orthogonal group $ O(7,\mathbb{R})=\{Q\in  GL(n,\mathbb{R}) | Q^{T} Q=Q Q^{T }=I \} $, $G_2 = \{g\in GL(n,\mathbb{R}) |g^* (\phi)=\phi \}$ a few properties of $G_2$ will be used here. The group $G_2$ acts irreducibly on $\mathbb{R}^7$  and acts transitively on the unit sphere $S^6 \subset \mathbb{R}^7 $ and preserve the metric and orientation for which the  basis $\{e_1 ,e_2 ,...,e_7 \}$ is an oriented orthonormal basis, the notation g and $ \langle .  ,.  \rangle $ will be used to refer the metric.
  \\The stabilizer subgroup of any non-zero vector in $\mathbb{R} ^7 $ is isomorphism to $SU(3)\subset SO(6)$, so that $S^6=G_2 / SU(3)$. Since SU(3) acts transitively on $S^5 \subset \mathbb{R}^6$, it follows that $G_2$ acts transitively on the set of orthonormal pairs of vectors in $\mathbb{R}^7$. However, $G_2$ does not act transitively on the set of orthonormal triples of vectors in $\mathbb{R}^7$ since it preserve the 3-form $\varphi=*\phi$.
  \\By definitions, $G_2$ is a compact lie group as a closed subgroup of the orthogonal group $O(7,\mathbb{R})$ and Lie algebra of $G_2$ is $\mathfrak {g}_2$ .

  \section{Nearly Half Flat SU(3)-structures}
\textbf{Definition.}Let $\mathfrak{g}$ a lie algebra then  $\mathfrak{g}$ is nilpotent if the lower central series terminates , i.e if  $\mathfrak{g}_n=0$  for some $ n \in \mathbb{N}.$ this means that 
$$[ X_1,X_2,[...[X_n,Y]...]]=ad_{X_{1}}ad_{X_{2}}...ad_{X_{n}}Y \in\mathfrak{g}_n$$ 
$$X_1,X_2,...,X_n,Y \in \mathfrak{g} $$
so that $ad_{X_{1}}ad_{X_{2}}...ad_{X_{n}}=0.$
Note that if $\mathfrak{gl}(n,\mathbb{K})$ is the set of $n\times{n}$ with the entries of $\mathbb{K}$, then the subgroup consisting of strictly upper triangular matrices is nilpotent lie algebras. 
\\A Lie derivative for a general differential form is likewise a contraction, taking into account the variant in $X:L_{x}w=\iota_{x}dw+d(\iota_{x}w)$ knowing as Cartan's Formula.\vspace{1cm}\\ 
\textbf{Definition.} For a lie algebra $\mathfrak{h}$ of dimension 2m, a closed two-form $\omega \in \Lambda^2 \mathfrak{h} ^{*} $ is called symplectic if it is non-degenerate i.e. in the case of $m=3$ $\omega ^3 \neq 0$. \vspace{4mm}\\

\textbf{Note1.}From the sixtieth classes of $G_2$-structures of in this paper we consider the class of torsion $W_1\oplus W_3\oplus W_{4}$ .\\
Gray and Harvella in [4] proved that there exist sixteen different class of almost Hermitian structures attending to the behavior of the covariant derivative of its $K\ddot{a}ler$ form. equivalently the different classes of SU(3)-structures can be defined in terms of form $\omega, \psi_{+}$ and $\psi_{-}$ in particular we are interested on Nearly Half-Flat which is defined on six-dimension lie algebra $\mathfrak{h}$.\\ 
An SU(3)-structure on a manifold of real dimension 6 consists of a Hermitian structure (g, J, $\omega$) and a $(0,3)-$form $\psi$; and SU(3) is the stabilizer of the transitive action of the group $G_2$ on $S_6$.\\

\textbf{Definition.} We call an SU(3)-structures $(\omega,\psi_+,\psi_{-} )$ is nearly half-flat on a 6-dimensional manifold when it satisfied the equation: 
$$d\psi_{-}=\dfrac{1}{2}(\omega\wedge\omega) $$
However six-dimensional nilpotent lie algebras admitting a nearly half-flat SU(3)-structures are not classified,Those which admitting a double half-flat are classified by Choissi and Swann in [1] and as double half-flat are in particular nearly half-flat this allows one to shows existence locally conformally cocalibrated $G_2$-structure.
\newpage

\section{Locally conformal cocalibrated $G_2$-structures}

Now we consider the L.C.CC $G_2$-structure which is introduce by a 4-form $\phi$ s.t. $d\phi=\theta \wedge \phi$ for a closed 1-form $\theta$, and $\phi \in \Lambda ^{4}V^{*}$ compatible with orientation of $\mathfrak{g}$ and the underlying metric g, consider the  $X \in\mathfrak{z(h)}$,  $L_{X}=0$ where $[X,Y]=L_{X} Y$  $X=\theta ^{\sharp}$ is isomorphism induced by the metric between the space and dual space, where:  
$$d\phi=\theta \wedge \phi, X=\theta ^{\sharp} \in \mathfrak{g} $$ s.t. $\phi= \sigma+\psi_{-} \wedge \theta$ , $\iota _{X}\omega^{2} =0$ and $\iota _{Y}\psi_{-}=0$ $\Rightarrow \phi=\theta \wedge \phi  \\\Rightarrow d\omega^{2} = (\omega^{2} -d\psi_{-})\wedge\alpha$ \vspace{15mm} \\observe that\vspace{5mm} $\theta(X/{\rvert X \rvert})=1$ for $X\in \mathfrak{z(h)} $ and $L_{X}=0$, $[X,Y]=L_{X} Y $ \\where $\rvert X \rvert = \sqrt{g(X,Y)}$ for $X^{\perp}= \{Y \in \mathfrak{g} \mid g(X,Y)=0 \}$ , \vspace{5mm} 
$$\theta=\theta _{1}e^1+\theta _{2}e^2+...+\theta _{6}e^6+\theta _{7}e^7$$ $L_{X}\phi
=\iota _{X}d\phi + d\iota _{X} \phi$, $L_{X}\phi=\phi + d(-\psi _{-})$ where $\psi _{-}=-\iota _{X} \phi$ where pair $(\omega,\psi_{-})$ is defining a neraly half-falat SU(3)-structure,
\\$de^i (e_l,e_k)=-e^i ([e_l,e_k])$,  $X=\theta ^{\sharp}$, $\theta(Y)=g(X,Y)=\theta ^{\sharp} for \forall Y\in \mathfrak{g} $ theta is dual form of vector, s.t. $\theta \in \Lambda ^{1} \mathfrak{g}^{*}$ our vector $X \in[\mathfrak{g},\mathfrak{g}]$ where $ \Lambda ^{1} \mathfrak{g} \longrightarrow \Lambda ^{1} \mathfrak{g}^{*} $, $X\in [\mathfrak{g},\mathfrak{g}]$ and $G_2$ acts irreducibly on V and hence on $\Lambda ^{1}V^{*}$
where there is a short exact sequence: $$0\longrightarrow \mathbb{R}X\longrightarrow \mathfrak{g} \longrightarrow \mathfrak{h}\longrightarrow0  $$  with the definition of this homomorphism $\pi^{*}:\Lambda ^{k} \mathfrak{g}^{*}\longrightarrow \Lambda ^{1} \mathfrak{h}^{*}$ it send $\pi^{*} \alpha$ for  one-form $\alpha$ s.t. $$\pi^{*} [X,Y]=[\pi^{*}X,\pi^{*}Y]$$\\
\newpage

\begin{prop}\textbf{.}
Consider $\phi$ be a 4-form Locally Conformal Cocalibrated $G_2$-structures i.e. $d\phi=\theta \wedge \phi$ ,$d\theta =0$ , $\theta \in \Lambda ^{1}\mathfrak{g}^{*}$ and the pair of $(\omega,\psi _{-})$ on  $\mathfrak{h}$, and by definition $\psi _{-}:= -\iota _{X}\phi$ , $\sigma:= \phi-\psi _{-} \wedge \theta$, $\sigma:= \omega\wedge\omega$ for $\psi _{-}\in \Lambda ^{3}\mathfrak{h}^{*}$, $\omega\in \Lambda ^{2}\mathfrak{h}^{*}$ ,then 
The pair $(\omega,\psi _{-})$ has a nearly half-flat SU(3)-structures which $d\psi _{-}=k \omega^{2}$ for $k \in \mathbb{R}-\{0\}$ with the torsion type $W^{\pm}_1\oplus W^{-}_2\oplus W_{3}$.\vspace{1cm} \\
\textbf{Proof.}\vspace{1cm}\\ (I) Contraction of $\psi _{-}$ and $\sigma:= \omega\wedge\omega$ 
\begin{itemize}
\item $\iota_{X}\psi_{-}=\iota_{X}(-\iota _{X}\phi)=-\iota _{X}\iota _{X}\phi=0$

\item $\iota _{X}\sigma=\iota _{X}\phi-\iota _{X}(\psi _{-}\wedge\theta)
              =-\psi _{-}-(\iota _{X}\psi _{-}\wedge\theta-\psi _{-}\wedge\iota _{X}\theta)=-\psi _{-}+\psi _{-}\wedge \iota _{X}\theta=-\psi _{-}+\psi _{-}=0$
 
\item $X=\theta ^{\sharp}$ and $\iota _{X}\theta=\theta(X)=1$ where $\theta(Y)=g(X,Y)$ $\forall Y\in\mathfrak{g}$\\ $\theta(X)=g(X,X)=|X|^2$\\
 
\end{itemize}(II) The derivative of 3-form $\psi _{-}$ and $\sigma:= \omega\wedge\omega$ are as follow, and it admits nearly half-flat SU(3)-structures. 
\begin{itemize}

\item $d\psi _{-}=d(-\iota _{X}\phi)=\iota _{X}(d\phi)=\iota _{X}(\theta\wedge\phi)=(\iota _{X}\theta)\wedge\phi-\theta\wedge(\iota _{X}\phi)=\theta(X)\phi-\theta\wedge(-\psi _{-})=\phi+\theta\wedge\psi _{-}=\sigma+\psi _{-}\wedge\theta+\theta\wedge\psi _{-}=\sigma$ $\Longrightarrow$ $d\psi _{-}=\sigma =\dfrac{1}{2}\omega^2$ $\Longrightarrow$ $k=\dfrac{1}{2}$ , $k \in \mathbb{R}-\{0\}$

\item $d\sigma=d(\phi-\psi _{-}\wedge\theta)=d\phi-d(\psi _{-}\wedge\theta)=d\phi-d\psi _{-}\wedge\theta+\psi _{-}\wedge d\theta=\theta\wedge\phi+\theta\wedge d\psi _{-}=\theta\wedge(\phi+d\psi _{-})=\theta\wedge(\phi+\sigma)$

\end{itemize}

\end{prop}\vspace{1cm} 

\begin{thm}\textbf{.}
Let $\mathfrak{h}$ be a six-dimensional nilpotent Lie algebra admitting a Nearly Half-Flat SU(3)-structures given by the pair $(\omega,\psi_{-})$,and a closed Lee form, $\theta$, such that $\theta(X)=1$,  $\forall X\in \mathfrak{z(h)}$,then the seven-dimensional Lie algebra $\mathfrak{g}=\mathfrak{h}\oplus \mathbb{R} X $ admits a Locally Conformal Cocalibrated $G_2$-structures.\vspace{1cm} \\

\textbf{Proof.} For a four-form $\phi=\frac{1}{2}\omega^2+\psi _{-}\wedge \theta$ , $\phi \in \Lambda ^{4} \mathfrak{g}^{*}$ Look at parts (I) and (II) of preposition 1.\vspace{1cm}

\end{thm}
\newpage 
\textbf{Corollary 1.} Those six-dimensional Lie algebras $\mathfrak{h}$ with double half-flat which are classified in \cite{swann} are admitting nearly half-flat SU(3)-structures for $X\in \mathfrak{z(h)}$, then $\mathfrak{g}=\mathfrak{h}\oplus \mathbb{R} X $ where $\phi \in \Lambda ^{4}\mathfrak{g}^{*}$ is Locally Conformally Cocalibrated $G_2$-structures, if $\mathfrak{g}=1A,1B,1C,2C,2B$ and $3A$.\vspace{1cm}\\
$$1A=(0,0,0,e^{12},e^{13},e^{23})$$
$$3A=(0,0,0,0,e^{12},e^{15+34})$$
$$1B=(0,0,e^{12},e^{13},e^{23},e^{14})$$
$$2B=(0,0,e^{12},e^{23},e^{14+35})$$
$$1C=(0,0,e^{12},e^{13},e^{23},e^{14+25})$$
$$2C=(0,0,e^{12},e^{13},e^{23},e^{14-25})$$
 \vspace{1cm}\\
 \textbf{Proof.} The set of double half-flat six-dimensional lie algebras are subset of  nearly half-flat six-dimensional lie algebras as any lie algebra $\mathfrak{h_{0}}$ with nearly half-flat SU(3)-structures has a equality in which $d\psi_{-}=\dfrac{1}{2}(\omega\wedge\omega)$, where in the case of double half-flat we have this equality in pulse : $d\psi_{+}=0$ .\\(section 4 of paper $[1]$) \vspace{1cm}\\
 This allows ones to conclude that there exist L.C.CC $G_2$-structures on seven dimensional decomposable nilpotent lie algebras of the form $\mathfrak{g}=\mathfrak{h}+\mathbb{R}e_{7}$ where $\mathfrak{h}$ admits double half-flat lie algebras, so we are able to state the existence result for six double half-flat SU(3)-structures.But these are first result for building L.C.CC $G_2$-structures, as up to now we didn't find a stable 4-form $\phi_{0}$ compatible with the orientation and metric g and it could be the next step on this paper.

\newpage\section{Associated Lie algebras}
 
 Giving a method to obtain new 7-dimensional lie algebras endowed with locally conformally co-calibrated (L.C.CC) $G_2$-structures. Starting from 6-dimensions lie algebras with the nearly half-flat SU(3)-structures and then describing some explicit examples of them.
 \\ A SU(3)-structure on a lie algebra $\mathfrak{h} $ of dimension six, consist in a triple $(g,J,\Psi)$ that $(g,J)$ are almost Hermitian structure on $\mathfrak{h} $ and $\Psi =\psi _{+}+i\psi _ {-}$ is a complex volume $(0,3)$-form and satisfying :
 $$ \dfrac{1}{6!}\omega ^6={(-1)}^3(\dfrac{i}{2})^3 \psi \wedge \bar{\psi}$$ 
 The $K\ddot{a}ler$ form associated to $(\omega,\psi _{-})$ such that describe a metric as $$g(X,Y)\omega ^{3}=-3 \iota _{X}\omega \wedge \iota _{Y}\psi _{-} \wedge \psi_{-} $$ with $X,Y \in \mathfrak{h}$ and $\iota _{X}$ denoting the contraction by X, when $(g,J,\psi _{-})$ is a nearly half-flat SU(3)-structures on lie algebras $\mathfrak{h}$ , we may choose  an orthonormal frame $\{e_1,e_2,...,e_6\}$  for almost complex structure s.t. $\psi_{-}=J \psi _{+}$
 with the $\{e^1,e^2,...,e^6\}$ an orthonormal dual basis and form $\omega$ and 3-form $\psi_{-}$ can be written as   
 $d\omega =0$ ,
 $\Psi = \psi_{+}+ i\psi _{-}$

 Now consider $\phi$ is a locally conformally co-calibrated $G_2$-structures and restricting the condition to decomposable cases s.t. $$\mathfrak{g} =\mathfrak{h}\oplus \mathbb{R}$$ again consider $\{e_1,e_2,...,e_6\}$ the basis for six-dimension vector space V and $\{e^1,e^2,...,e^6\}$ basis for dual space $V^{*}$ ,with one form $\theta \in \varLambda ^{1} V^{*}$ and pair of $(\omega,\psi_{-})$ nearly half-flat SU(3)-structures where $\omega \in \varLambda ^{2} V^{*}$ and $\psi_{-} \in \varLambda ^{3} V^{*}$ on six-dimension lie algebras $\mathfrak{h}$, for $20$ possible different ways to have $3$-dimension basis as $ \dfrac{6!}{3!3!}=20 $
 
 $$\begin{array}{c  c c c c} 
 \setlength{\tabcolsep}{1.0 cm}
 	e^{123}  & e^{124} & e^{125} & e^{126} & e^{134}
 	\\e^{135} & e^{136} & e^{145} & e^{146} & e^{156}
 	\\e^{234} & e^{235} & e^{236} & e^{245} & e^{246}
 	\\e^{256} & e^{345} & e^{346} & e^{356} & e^{456}
 \end{array} $$ 
 Where $e^{ijk}=e^{i} \wedge e^{j} \wedge e^{k}$ for $i,j,k \in \{1,2,...,6\}$, are ordered as $ i<j<k $\vspace{1cm} $$\psi_{-}:=\sum_{i<j<k}\psi_{ijk}e^{ijk}$$.\newpage
 
 Is a generic 3-form.
 For four-form $\varphi \in \Lambda^{4}V^*$  where $V^*$  is dual space of V in seven-dimension with the basis $\{e^1,e^2,...,e^6,,e^7 \}$ and $\{e_1,e_2,...,e_6,e_7 \}$ respectively, so all of the possibilities for different basis would be $ \dfrac{7!}{3!4!}=35 $ 
 
 $$\begin{array}{c  c c c c} 
 
 e^{1234}  & e^{1235} & e^{1236} & e^{1237} & e^{1245}
 \\e^{1246} & e^{1247} & e^{1256} & e^{1257} & e^{1267}
 \\e^{1345} & e^{1346} & e^{1347} & e^{1356} & e^{1357}
 \\e^{1367} & e^{1456} & e^{1457} & e^{1467} & e^{1567}
 \\e^{2345} & e^{2346} & e^{2347} & e^{2356} & e^{2357}
 \\e^{2367} & e^{2456} & e^{2457} & e^{2467} & e^{2567}
 \\e^{3456} & e^{3457} & e^{3467} & e^{3567} & e^{4567}
 \end{array} $$ \vspace{5mm}\\

 \textbf{\underline{Two examples of L.C.CC $G_2$-structures}:}\\

 First consider the lie algebra $1A=(0,0,0,e^{12},e^{13},e^{23})$ and one form $$\theta:= a_1e_1+a_2e_2+a_3e_3+a_4e_4+a_5e_5+a_6e_6+a_7e_7$$ where $\theta \in \Lambda^{1}V^*$
 then the 4-form $\Phi$ when the generic form is : $$\phi:=\sum_{i<j<k<l}\phi_{ijkl}e^{ijkl}$$
 and $\phi \in  \Lambda^{4}V^* $\vspace{5mm}
 
$\phi= \phi_{3567}e^{3567}+\phi_{1467}e^{1467}+\phi_{1567}e^{1567}+\phi_{2345}+\phi_{2347}e^{2347}+\phi_{2357}e^{2357}+\phi_{2367}e^{2367}+\phi_{2467}e^{2467}+\phi_{2567}e^{2567}+\phi_{1257}e^{1257}+\phi_{1267}+\phi_{1347}e^{1347}+\phi_{1357}e^{1357}+\phi_{1367}e^{1367}+\phi_{1457}e^{1457}+\phi_{1234}e^{1234}+\phi_{1235}e^{1235}+\phi_{1236}e^{1236}+\phi_{1237}e^{1247}+a_{7}\phi_{2345}e^{4567}+(-a_{7}\phi_{1235}+\phi_{1567})e^{3457}+(-a_{7}\phi_{1236}+\phi_{2567})e^{3467}+(-a_{7}\phi_{1234}+\phi_{1467})e^{2457}+\phi_{2345}e^{1256}-\phi_{2345}e^{1346} $\\

   This is a 4-form in seven dimension with the real lie algebra 1A, this form accepted locally conformal cocalibrated (L.C.CC) $G_2$-structure.\vspace{5mm}\\
   Now consider third , $1B=(0,0,e^{12},e^{13},e^{14})$ and one form $\theta$ as previous again: $$\theta:= a_1e_1+a_2e_2+a_3e_3+a_4e_4+a_5e_5+a_6e_6+a_7e_7$$\\
   
   $\phi= (a_7\phi_{1236}+\phi_{1567})e^{2467}+\phi_{2567}e^{2567}+(a_7\phi_{1245}+\phi_{2567})^{3457}+(a_7\phi_{1246})e^{3467}+\phi_{1457}e^{1457}+\phi_{1467}e^{1467}+\phi_{1567}e^{1567}\phi_{2567}e^{2567}+\phi_{1257}e^{1257}+\phi_{1267}+\phi_{1347}e^{1347}+\phi_{1357}e^{1357}+\phi_{1367}e^{1367}+\phi_{1457}e^{1457}+\phi_{1234}e^{1234}+\phi_{1235}e^{1235}+\phi_{1236}e^{1236}+\phi_{1237}e^{1247}+a_{7}\phi_{2345}e^{4567}+(-a_{7}\phi_{1235}+\phi_{1567})e^{3457}+(-a_{7}\phi_{1236}+\phi_{2567})e^{3467} $
\newpage
\textbf{Acknowlegment}:
 I would like to thanks Alberto Raffero and Leonardo Bagaglini and very special thanks to Anna Maria Fino for giving the problem and her suggestions, also I am thankful for Sigbj\o{}rn Hervik for his interest on my work. 
 \vspace{1cm}\\
\textbf{Contact.} Dipartamento di matematica,Via Carlo Alberto 10,\\ Universit\`{a} degli studenti di Torino, Italy.
E-mail: azohrabi@unito.it

\end{document}